\documentclass[12pt]{amsart}
\usepackage{amsmath,amssymb,amsthm,amscd}

\title[Pattern equivariant functions and deformations]{Pattern
  equivariant functions, deformations and equivalence of 
tiling spaces}
\author{Johannes Kellendonk}
\address{Université de Lyon,
Université Lyon1,
CNRS, UMR 5208 Institut Camille Jordan,
Batiment du Doyen Jean Braconnier,
43, blvd du 11 novembre 1918,
F - 69622  Villeurbanne Cedex,
France}                                         

\date{\today}

\newcommand{\N}{\mathbb{N}}

\newcommand{\R}{\mathbb{R}}
\newcommand{\C}{\mathbb{C}}

\newcommand{\rn}{\R^n}
\newcommand{\rem}{\R^m}
\newcommand{\rphi}{\Phi}
\newcommand{\tphi}{\Phi_\varphi}
\newtheorem{thm}{Theorem}[section]
\newtheorem{lemma}[thm]{Lemma}
\newtheorem{cor}[thm]{Corollary}
\newtheorem{prop}[thm]{Proposition}
\theoremstyle{definition}
\newtheorem{defn}[thm]{Definition}
\newtheorem{hypo}[thm]{Hypothesis}
\theoremstyle{remark}
\newtheorem{remark}[thm]{Remark}

\numberwithin{equation}{section}

\newcommand{\flc}{finite local complexity}
\newcommand{\Pp}{\mathcal P}

\newcommand{\cp}{\Gamma}
\newcommand{\Om}{\Omega}

\newcommand{\eP}{$P$-equivariant}
\newcommand{\sP}{strongly $P$-equivariant}
\newcommand{\wP}{weakly $P$-equivariant}
\newcommand{\dash}{\!-\!}

\newcommand{\fid}{{\hbox{\rm\footnotesize id}}}
\newcommand{\orb}{{\hbox{\rm orb}}}

\begin{document}

\date{\today}

\maketitle

\newcommand{\supp}{\mbox{\rm supp}}
\newcommand{\id}{\mbox{\rm id}}



\begin{abstract} We reinvestigate 
the theory of deformations of tilings using
$P$-equivariant cohomology. In particular we relate the 
notion of asymptotically negligible
shape functions introduced by Clark and Sadun to weakly $P$-equivariant
forms. We then investigate more closely
the relation between deformations of patterns and homeomorphism
or topological conjugacy of pattern spaces.
\end{abstract}

\section{Introduction}
The study of aperiodic systems in physics or geometry
has led to the definition of cohomology groups associated
with aperiodic tilings or point sets of $\rn$ (we use here the word pattern
to mean either of them)
In physics some elements of these groups are
related (via $K$-theory and cyclic cohomology) 
to topologically quantized transport
properties, see \cite{KellendonkRichard} for a recent overview. 
In geometry Sadun, Williams and 
Clark have given an interpretation of the
cohomology group (with values in $\rn$) 
in terms of deformation theory of tilings \cite{SadunWilliams03,ClarkSadun06}.
In short, an (admissible) 
$1$-cocycle defines a deformation of a tiling by redefining
the shape of its tiles.
If the cocycle is a coboundary
then the deformed tiling is locally derivable from the
original one.   
A deformation alters the properties of the dynamical system
associated with the tiling,
except if the new tiling is topologically conjugate to the old one, a notion
which is, however, strictly weaker than being mutually locally derivable.   
To capture also topological conjugacy in cohomological terms, Clark and Sadun 
introduced the concept of asymptotically negligible cocycles and
showed that such cocycles yield deformations which are topologically
conjugate.

Clark and Sadun used a formulation of tiling cohomology which is based
on the Anderson-Putnam-G\"ahler construction. This construction
furnishes a system of finite CW-complexes so that one can make use 
of their cellular cohomology. 
Our aim here is to provide a formulation of deformation theory in terms
of pattern equivariant functions making use of de Rham cohomology.
From this point of view a deformation of a pattern is given by a
pattern equivariant $1$-form.
This allows us to formulate the somewhat ad hoc notion of
asymptotically negligible cocycles in terms of natural analytic
properties of pattern equivariant functions (Theorem~\ref{thm-neg}). 
In fact, as was
introduced in \cite{KellendonkPutnam06}, there are naturally strongly
and weakly pattern equivariant forms on $\rn$ and we find that the
theory of deformation modulo topological conjugacy of a pattern $P$
is related to the mixed cohomology group,  
$$ Z^1_{s\dash P}(\R^n,\R^n) / 
B^1_{w\dash P}(\R^n,\R^n)\cap Z^1_{s\dash P}(\R^n,\R^n),$$
where $ Z^1_{s\dash P}(\R^n,\R^n)$ denotes the $\rn$-valued strongly
$P$-equivariant closed $1$-forms and $B^1_{w\dash P}(\R^n,\R^n)$ the 
$\rn$-valued weakly $P$-equivariant exact $1$-forms.
A study of the mixed cohomology group for certain classes of tilings
is under investigation, for substitution tilings the results of
\cite{ClarkSadun06} point out that this is worthwhile.

As has already been observed by Sadun and Williams, deformations give
rise to homeomorphisms between tiling spaces. We review this result in
the framework of $P$-equivariant functions. We find that
homeomorphisms coming from deformations preserve the canonical
transversals and show that this property 
characterises such homeomorphisms: any homeomorphism between two pattern
spaces preserving the canonical transversals and identifying the two
patterns comes from a deformation (Cor.~\ref{cor-def}).

We further investigate deformations which give (topologically)
conjugate patterns. Our Theorem~\ref{thm-exact} 
can be seen as the analog of Clark and Sadun's result about
asymptotically negligible cocycles: deformations differing from the
identity by a \wP\ $1$-form
lead to conjugate patterns. More presisely we find that they give rise
to a second homeomorphism which also identifies $P$ with its deformed
pattern, but does not preserve the canonical transversal. Instead it
commutes with
the $\rn$-actions. Under an additional assumption on the deformation,
which we call boundedness, a converse can be obtained: bounded
deformations which yield (pointed) conjugate patterns differ from the identity
by a \wP\ $1$-form (Theorem~\ref{thm-exact-converse}).

Finally we present a detailed analysis of the question of
invertibility of deformations 
(which somehow is hidden 
in the notion of admissibility in \cite{ClarkSadun06}). 
We show that there exists a neighbourhood of the identity deformation
which contains only invertible deformations.

\section{Preliminaries}
\subsection{Pattern spaces and their dynamical systems}\label{sec-2.1}
The objects we are interested in and which we want somehow to compare are 
closed subsets of euclidean space, mostly of \flc. These might either
be uniformly discrete or they might be tilings; for shortness we call
them patterns.

Let $P,Q$ be two closed subsets of $\rn$. Define their distance to be
$$D(P,Q)= \inf \left\{\epsilon>0 \left| d_H(B_{\frac{1}{\epsilon}}[P]\cup
    \partial B(0,\textstyle{\frac{1}{\epsilon}}), 
B_{\frac{1}{\epsilon}}[Q]\cup \partial B(0,\textstyle{\frac{1}{\epsilon}}))\leq
\epsilon\right.\right\}. $$ 
Here $B(x,r)$ is the (open) ball of radius $r$ around $x$,  
$\partial B(x,r)$ denotes its boundary,
$$B_r[P] = B(0,r)\cap P,$$
and $d_H$ is the Hausdorff distance. We call $B_r[P]$, $x\in\rn$,
the $r$-patch of $P$ around $0$.

We will mostly look at discrete subsets of $\rn$ but for comparison
with \cite{ClarkSadun06} we also look at tilings.
There are different
ways to associate a discrete point set to a tiling $T$. A
possibility would be to introduce a point in the interieur of each
tile, in such a
way that translationally congruent tiles have that point at the same
position. The set of all such points is called the set 
of punctures of the tiling. We denote by $T^{pt}$. Another possibility
which turns out usefull for polyhedral tilings is to look at the set
of their vertices, denoted $T^0$. The two different ways have their
own advantages. For instance, there is a natural
bijection between the tiles of $T$ and the points of $T^{pt}$.

$P$ is uniformly  discrete if its points have a minimal distance, i.e.\
if
$$ r_{min}(P) = \inf\{|x-y| : x\neq y \in P\} $$
is strictly positive. If $h\in\rn$ and 
$|h|<\frac{r_{min}}{2}$ then $D(P,P-h)\leq |h|$. 

$P$ is relatively dense if it does not have arbitrarily large holes,
i.e.\ there exists $r$ such that for all $x\in\rn$, $B_r[P-x]$
contains at least one point. 
A uniformly discrete and relatively dense set is called a Delone set.

For a closed subset $P\subset\rn$ the completion of $\orb(P):=\{P-x|x\in\rn\}$,
the $\rn$-orbit of
$P$, is called the continuous hull of $P$, written $\Omega_P$,
$$\Omega_P := \overline{\orb(P)}^D.$$
$\Omega_P$ is a compact space which 
carries a continuous action of $\rn$ induced by the translation action
on $P$. We denote the action of $x$ by $\omega\mapsto
\omega-x$. This defines the (continuous) dynamical system
$(\Omega_P,\rn)$ of $P$. 

If $P$ is uniformly discrete we consider also its discrete hull
$$\Xi_P := \overline{\{P-p|p\in P\}}^D.$$
It is a transversal for the $\rn$ action and also referred to as the
canonical transversal of $P$. 
If $T$ is a tiling one may define a transversal and an $r$-discrete groupoid
either with $T^0$ or with $T^{pt}$. These would differ for the two
different choices, but our results below apply to both, provided they
are uniformly discrete.

$P$ (or $T$) has \flc\ (w.r.t.\ the translation group) 
if for any given $r$ there are only finitely many distinct
$r$-patches $B_r[P-x]$ with $x\in P$ (or $x\in T^{pt}$ or $x\in T^0$).

Finite local complexity has strong implications on the structure of
the hull. First, all elements of the hull of a discrete pattern $P$
with \flc\ can be viewed as uniformly discrete subsets of $\rn$ (including,
if $P$ is not relatively dense, the empty set). Second, such a hull
carries an equivalent
but somewhat simpler metric, namely
$$D_t(Q,Q') =  \inf_{r>0}\left\{\left. r^{-1}\right|
B_r[Q] \cong_{r^{-1}} B_r[Q'] \right\}, $$
where we write $B_r[Q] \cong_{r^{-1}} B_r(Q')$ for
$\exists x,x'\in
B(0,\frac{1}{2r}): B_r[Q-x] = B_r[Q'-x']$.
Restricted to the canonical transversal the metric is even equivalent to
$$D_0(Q,Q') =  \inf_{r>0}\left\{\left. r^{-1}\right|
B_r[Q] = B_r[Q'] \right\}. $$
Third, the canonical transversal $\Xi_P$ is totally disconnected.
Finally, at least if $P$ is also relatively dense, 
the hull is a foliated space whose
leaves are locally homeomorphic to $\rn$. To see this let 
$U_\epsilon(Q)$ be the open $\epsilon$-ball around $Q\in \Omega_P$ w.r.t.\
the metric $D_t$ and $T_{Q,\epsilon,V}=\{Q'-v| B_{\frac{1}{\epsilon}}[Q'] = 
B_{\frac{1}{\epsilon}}[Q],\:v\in V\}$ where $V\subset \rn$. 
Then for any small
enough $\epsilon$ one finds $\epsilon_1,\epsilon_2>0$ such that 
$U_{\epsilon_1}(Q)\subset
T_{Q,\epsilon,B(0,\epsilon)} \subset U_{\epsilon_2}(Q)$.
If $P$ is relatively dense then translates $U_\epsilon(Q)-x$ with
$Q\in\Xi_P$ and $|x|<r_{max}(P)$
generate the topology and hence the family of sets
$T_{Q,\epsilon,B(x,\epsilon)}$, $Q\in\Xi_P$, $|x|<r_{max}(P)$,
$\epsilon<\frac{r_{min}(P)}{2}$ generates the topology as well.
Now $T_{Q,\epsilon,B(x,\epsilon)}$ is the 
homeomorphic image of $T_{Q,\epsilon,\{0\}}\times B(x,\epsilon)$ under the map
$(Q',v)\mapsto Q'-v$. We thus have charts which are
Cartesian products of totally
disconnected sets (clopen subsets of $\Xi_P$) with open balls of $\rn$. 
These are the foliation charts.
This has two consequences which will be
important further down: First, the path-connected components of $\Omega_P$
are the orbits under the translation action, and second, the pre-image under
the map $\rn\to\orb(P)$: $x\mapsto P-x$
of the path-connected
component of $P$ in $\orb(P)\cap U_{\epsilon}(P)$ lies in
$B(0,\epsilon)$.
\subsection{Notions of equivalence}
There are a variety of equivalence relations between patterns. Some of them 
are expressed directly between the patterns, as for instance, 
mutual local derivability:
\begin{defn}[ \cite{BJS}] Let $P$ and $Q$ be subsets of $\rn$. 
$Q$ is locally derivable from $P$ if 
there exists an $R>0$ such that for all $x,y\in\rn$:
\begin{equation}\nonumber
B_{R}[P-x] = B_{R}[P-y] \Longrightarrow 
\{0\}\cap (Q-x) = \{0\}\cap (Q-y)\, .
\end{equation}
We say that $R$ is a derivability range of the local derivation.
\end{defn}
This definition applies to tilings if we regard them as subsets of
$\rn$ being given by the boundary points of their tiles.
 
Equivalently, $Q$ is locally derivable from $P$ if and only if
for all $r>0$ exists $r'$ such that for all $x,y\in\rn$:
\begin{equation}\nonumber
B_{r'}[P-x] = B_{r'}[P-y] \Longrightarrow 
B_r[Q - x] = B_r[Q - y]\, .
\end{equation}
This has the following interpretation which explains the name: One can
construct the $r$-patch of $Q$ around $x$ from the 
$r'$-patch of $P$ around $x$.

Other notions of equivalence refer to the spaces and dynamical systems
defined by the patterns. 
For instance, we could regard patterns as equivalent if they have 
\begin{enumerate}
\item\label{re2} 
homeomorphic hulls, or
\item\label{re3} topological conjugate (continuous) dynamical systems.
\end{enumerate}
Relation \ref{re3} implies relation \ref{re2} 
(trivally). In fact, we will find a stronger relation associated with 
a deformation, namely an invertible deformation $P'$ of $P$ (we
explain the notion further down) leads to a homeomorphism between the
hulls which preserves the canonical transversal and identifies $P$ with $P'$.
Likewise we find a relation stronger than \ref{re3}, namely we say
that $P'$ is pointed conjugate to $P$ if there is a
topological conjugacy between the (continuous) dynamical systems which
maps $P$ to $P'$. We will also say that 
$P'$ is pointed semi-conjugate to $P$ if there is a
continuous surjection $\Omega_P\to\Omega_{P'}$ which commutes with the
action of $\rn$ and maps $P$ to $P'$.
Note that pointed conjugacy implies homeomorphic hulls
but the conjugacy will not, in general,
preserve the canonical transversal.

\subsection{$P$-equivariant functions}
Topological spaces can equivalently be decribed by commutative
$C^*$-algebras. Here these algebras consist of \eP\
functions. For that we need to assume that our subset $P$ has \flc.

Roughly speaking, 
a function $f$ on $\rn$ is $P$-equivariant with range $R$ if it
satisfies the following property:
if the patterns in $P$ surrounding two points $x$ and $y$ 
are equal in a ball of radius $R$
(after translating by $-x$ and $-y$, respectively), then $f$ must take the 
same values at $x$ and  $y$. More precisely, we have the following.
\begin{defn}[\cite{Kellendonk03}]
Let $f$ be a function defined on $\rn$. We say that $f$ is
$P$-equivariant with range $R$ if for $x,y\in \rn$:
$$
B_R[P - x]  = B_R[P - y] \Longrightarrow f(x) = f(y).
$$ 
$f$ is \sP\ if it is $P$-equivariant with finite range.
We also call a function $f$ defined on a subset $Q \subset \rn$
which is locally derivable from $P$ 
\sP\ if it satisfies the above conditions for $x,y\in Q$.
\end{defn}
Note that $Q$ is locally derivable from $P$ if and only if the indicator
function on $Q$ is \sP. 

A priori our definition does not require a \sP\ function to satisfy
any further regularity conditions. By adding such conditions we arrive
naturally at the concept of \wP\ functions. Let $f:\rn\to\rem$ and
$$ s_k(f) = \sup\{|D^\alpha f(x)|: x\in\rn, |\alpha|\leq k\},$$
$\alpha=(\alpha_1,\cdots,\alpha_n)\in \N^n$,
where
$D^\alpha f = \frac{\partial^{\alpha_1}}{\partial^{\alpha_1}x_1}\cdots
\frac{\partial^{\alpha_n}}{\partial^{\alpha_n}x_n} f$, and 
$|\cdot|$ denotes any norm on $\R^m$ but $|\alpha|=\sum_i^n\alpha_i$. 
\begin{defn}
We denote by $C^k_{s\dash P}(\rn,\R^m)$, $k\in\N\cup\{\infty\}$ 
the space of \sP\ $C^k$-functions over $\rn$ with values in $\R^m$. 
We say that a $C^k$-function $f$ is \wP\ if, for all $\epsilon>0$ there
exists a \sP\ $C^k$-function $f_\epsilon$ such that
$s_k(f-f_\epsilon)<\epsilon$, if $k$ is finite, and
$\forall l\in \N$: $s_l(f-f_\epsilon)<\epsilon$, if $k=\infty$.
The space of \wP\ $C^k$-functions  over $\rn$ with values in $\R^m$
is denoted $C^k_{w\dash P}(\rn,\R^m)$.
\end{defn} 
Since \sP\ functions are bounded, $C^k_{w\dash
  P}(\rn,\R^m)$ is the closure of $C^k_{s\dash P}(\rn,\R^m)$
in the space of functions $f:\rn\to\R^m$ which have bounded continuous
derivatives up to order $k$
w.r.t.\ the topology induced by $s_k$,
if $k$ is finite, or by all $s_l$, $l\in \N$, if $k=\infty$.
Smooth strongly or weakly $P$-equivariant forms can then be identified
with elements of $C^\infty_{s\dash P}(\rn,\Lambda{\rn}^*)$ or
$C^\infty_{w\dash P}(\rn,\Lambda{\rn}^*)$, respectively.
As usual, $df = \sum_i \frac{\partial}{\partial x_i}f dx_i$, and we
view $\{dx_i|i=1,\cdots,n\}$ as a base for ${\rn}^*$. 

Complex-valued continuous \wP\ functions\footnote{these were also
simply called $P$-equivariant functions over $\rn$} 
form a $C^*$-algebra.
The map $f\mapsto f_P$, $f_P(x):=f(P-x)$
induces an algebra isomorphism between $C(\Omega_P)$ and 
$C^0_{w\dash P}(\rn,\C)$.
The compactness of $\Omega_P$ is reflected in the following lemma
whose simple proof is left to the reader: 
\begin{lemma}\label{lem-PE}
Let $P$ be a uniformly discrete set with \flc. 
If $f\in C^0_{w\dash P}(\rn,\R^m)$ then its image is compact.
If even $f\in C^0_{s\dash P}(\rn,\R^m)$ then its restriction to $P$ 
takes only finitely many values.
\end{lemma}
We now present some general results about \eP\ functions.
\begin{lemma}\label{lem-trivial}
Let $P$ be a uniformly discrete set of \flc.
Let $\varphi$ and $\psi$ be two differentiable functions
which coincide on $P$. Suppose that $d\varphi$ and
$d\psi$ are \sP. Then $\varphi-\psi$ is \sP.
\end{lemma}
\begin{proof} Let $\eta=\varphi-\psi$. Hence $\eta$ vanishes on $P$ and
$d\eta$ is \eP, say with range
$R$. Let $r=R+r_{min}$. Let $B_r[P-x]=B_r[P-y]$. Choose 
$p\in P$ such that $p-x\in B_{r_{min}}[P-x]$. Let $q=p-x+y$.
Then 
$q-y=p-x\in B_{r_{min}}[P-y]$ and, by $P$-equivariance of $d\eta$,
$$\eta(x)-\eta(p)=\int_p^x d\eta = \int_{q}^y d\eta=\eta(y)-\eta(q).$$
Since $\eta(p) = \eta(q)$,
$\eta$ is \sP.\end{proof}

Let $f:P\to\rn$ be a function.
For $h=P-P:=\{p-q |p,q\in P\}$, we define $\Delta_h f:P\cap(P-h)\to\rn$ by
$$\Delta_h f (x) = f(x+h)-f(x).$$  
\begin{lemma}\label{lem-extend}
Let $\psi:P\to\rem$ 
be a function on a uniformly discrete set $P$ with \flc.
\begin{enumerate}
\item \label{item-1}
If $\psi$ is \sP\ there exists a smooth
\sP\ function $\varphi:\rn\to\rem$ which restricts on $P$ to
$\psi$.
\item \label{item-2}
If, for all $\epsilon>0$ there exists a \sP\ function
$\psi_\epsilon:P\to\rem $ such that for all $p\in P$:
$|\psi(p)-\psi_\epsilon(p)|<\epsilon$, then
there exists a smooth
\wP\ function $\varphi:\rn\to\rem$ which restricts on $P$ to $\psi$.
\item \label{item-3}
Suppose that $P$ is moreover relatively dense. 
If for all $h\in P-P$, the function $\Delta_h \psi:P\cap(P-h)\to\rn$
is \sP\ then
there exists a smooth function $\varphi:\rn\to\rem$ which restricts on
$P$ to $\psi$ and has \sP\ differential $d\varphi$.
\end{enumerate}
\end{lemma}   
\begin{proof}
Let $\rho:\rn\to\R$ be a 
smooth positive function with 
support contained in $B(0,r)$, $r\leq \frac{r_{min}}{2}$. 
For the first two parts of the lemma we consider
$$\varphi(x)=\sum_{p\in P}\rho(x-p)\psi(p)$$
and suppose that $\rho(0)=1$ so that $\varphi(p)=\psi(p)$ for
$p\in P$. 
If $\psi$ is \eP\ with range
$R$ then $\varphi$ is 
smooth and \eP\ with range $R+\frac{r_{min}}{2}$. 

Suppose that $\psi$ satisfies the conditions of the second
part. Then,  by the first part, 
$\varphi_\epsilon(x)=\sum_{p\in P}\rho(x-p)\psi_\epsilon(p)$
defines a smooth \sP\ function. Now 
$$|D^\alpha(\varphi(x)-\varphi_\epsilon(x))|\leq 
\|\psi-\psi_\epsilon \|_\infty 
|\sum_{p\in p}D^\alpha\rho(x-p)|\leq 
\epsilon \|D^\alpha\rho\|_\infty$$
and hence $\varphi_\epsilon$ tends to $\varphi$ in the Fr\'echet
topology when $\epsilon$ tends to $0$.

For the third part consider the Voronoi domains $V_p$ for
$p\in P$. Let $A$ be an upper bound for their diameter, it is finite
by the relative density of $P$. 
Define for all $x$ in the interior $\mbox{int}(V_p)$ of
the domain 
\begin{equation}\label{eq-ext}
\tilde\phi(x) := \psi(p).
\end{equation}
Then $\tilde\phi$ is defined almost everywhere and
$$ \varphi:=\rho*\tilde\phi$$
is a smooth function where we suppose that $\int\rho = 1$ so that $\varphi$
coincides with $\psi$ on $P$
provided $r$ is small enough. 
Now suppose that for all $h=P-P$ with $|h|\leq A$, 
$\Delta_h\psi$ is \eP\ with
range $R>3A$, and hence, if $B_R(P - p)  = B_R(P - q)$ for some $p,q\in
P\cap P-h$ then, for almost all $x\in B(0,R)$,
$\tilde\phi(x+p+h)-\tilde\phi(x+p)= \tilde\phi(x+h+q))-\tilde\phi(x+q)$. 
This implies that the equation 
$\tilde\phi(x+p)-\psi(p)= \tilde
\phi(x+q))-\psi(q)$, which by (\ref{eq-ext}) holds for $x\in V_p-p$,
extends to $x\in B(0,2A)$. It follows that
$\tilde\phi(x+p+h)-\psi(p)= \tilde 
\phi(x+h+q)-\psi(q)$ for $x\in B(0,2A)$ 
and hence
$$(d\rho * \tilde\phi)(x+p+h) = (d\rho * (\tilde\phi -
\psi(q)+ \psi(p))(x+q+h)=
(d\rho * \tilde\phi )(x+q+h)$$
where we have used that $d\rho * (\psi(q)-\psi(p))=0$.
Since $d\varphi(x) = (d\rho * \tilde\phi)(x)$ the statement follows. 
\end{proof}
We finally provide a condition under which a differentiable function
with \sP\ differential is \wP.
\begin{prop} \label{prop-set}
Let $P\in\rn$ be a Delone set of \flc\ and 
$\varphi:\rn\to\R$ be a $C^k$-function, $1\leq k\in\N\cup\{\infty\}$, with \sP\
differential. Then $\varphi$ is \wP\ if and only if
for all $\epsilon$ exists a finite partition of $P$ into subsets 
$P_1,\cdots,P_N$
which are locally derivable from $P$ and such that
$\forall i\in\{1,\cdots,N\}$ $\forall p,q\in P_i$: 
$|\varphi(p)-\varphi(q)|<\epsilon$.
\end{prop}
\begin{proof} Suppose that  $\varphi$ is \wP.
Choose $\epsilon$ and a \sP\ function
$\varphi_\epsilon$ such that
$\|\varphi-\varphi_\epsilon\|_\infty<\epsilon$. In particular
$\varphi_\epsilon (P)$ is a finite set, lets say $\{x_1,\cdots,x_N\}$.
Define $P_i=\varphi^{-1}_\epsilon(x_i)\cap P$ which is certainly
locally derivable from $P$. The $P_i$ partition $P$ and by
construction for all $p,q\in P_i$ we have
$|\varphi(p)-\varphi(q)|<2\epsilon$. 

We prove the converse:
Given $\delta$ choose a finite partition 
$P_1,\cdots,P_N$ of sets which are locally derivable from $P$ and such that
$\forall i\in\{1,\cdots,N\}$ $\forall p,q\in P_i$: 
$|\varphi(p)-\varphi(q)|<\delta$.
Note that this and the strong $P$-equivariance of $d\varphi$ 
implies that $\varphi$ has bounded derivatives $D^\alpha\varphi$,
$|\alpha|\leq k$.
Choose a point $p_i$ in each $P_i$.
Consider the Voronoi domains $V_p$ of $p$ in $P$. 
Define for all $x$ in the interior of $V_p$ for
$p\in P_i$ by
$\sigma_\delta(x)=\varphi(p)-\varphi(p_i)$. 
Then $\sigma_\delta$ is defined almost everywhere
and constant on the interior of the Voronoi domains. 
Since the $P_i$ are locally derivable from $P$, 
$\varphi-\sigma_\delta$ is \sP. It is also
$\delta$-close to $\varphi$ where defined and
the aim is to smoothen it out.

Let $\rho:\rn\to\R$ be a smooth positive
function with support contained in $B(0,1)$
and $\int\rho = 1$. Let $\rho_\delta(x) = \delta^{-n}\rho(\frac{x}{\delta})$.
Define $\varphi_{\delta_1,\delta_2} =
\rho_{\delta_1}*(\varphi-\sigma_{\delta_2})$ which is by construction
a smooth \sP\ function. 
We claim that $\forall \epsilon>0$ and $l\leq k$
$\exists  \delta_1,\delta_2$ such that 
$s_l(|\varphi-\varphi_{\delta_1,\delta_2}|) < \epsilon$. This then
proves that we can approximate $\varphi$ in the relevant topology by
\sP\ functions.

To proof the claim we fix $\epsilon$ and $l\leq k$. We have
$$ \|\rho_{\delta_1}*\sigma_{\delta_2}\|_\alpha\leq
\|\sigma_{\delta_2}\|_\infty\|D^\alpha \rho_{\delta_1}\|_1 \leq 
\delta_2\delta_1^{-|\alpha|}\|D^\alpha\rho\|_1.$$
Hence
$$ s_l(|\varphi-\varphi_{\delta_1,\delta_2}|) \leq 
s_l(|\varphi-\rho_{\delta_1}*\varphi|)+\sum_{|\alpha|\leq l}
\delta_2\delta_1^{-|\alpha|}\|D^\alpha\rho\|_1.$$
Since $D^\alpha\varphi$ is bounded continuous for $|\alpha|\leq k$
we can find
$\delta_1$ such that $s_l(|\varphi-\rho_{\delta_1}*\varphi|)<
\frac{\epsilon}{2}$ and then
$\delta_2$ 
such that $\sum_{|\alpha|\leq l}\delta_2\delta_1^{-|\alpha|}\|D^\alpha\rho\|_1
<\frac{\epsilon}{2}$. 
\end{proof}

\subsection{Cohomology groups for patterns and tilings}
Let $\Pp$ be a uniformly discrete 
subset or a tiling, of finite local complextity. 
Associated with $\Pp$ is a cohomology group $H(\Pp,A)$, the pattern cohomology,
with coefficients in an abelian group $A$. It can be described in a
variety of equivalent ways:

\subsubsection{Cech-cohomology}
$H(\Pp,A)$ is the Cech-cohomology $\check H(\Om_\Pp,A)$ of the
continuous hull. 

\subsubsection{Cellular cohomology}\label{sec-3.1}
Let $T$ be a polyhedral tiling, that is a tiling whose tiles are
polyhedra and match face to face. We can view $T$ as an infinite
CW-complex for $\rn$ and suppose to have chosen an orientation of its
cells such that translational congruent cells have the same orientation,
We also assume that $T$ has finite local complexity.
The Anderson-Putnam-G\"ahler construction expresses $\Om_\Pp$ as
  an inverse limit of CW-complexes. We recall the construction:
The $k$-neighbourhood of a tile $t$ in a polyhedral tiling $T$ is the
patch of all tiles of $T$ which meet the $\!k-\!1$-neighbourhood of $t$.
Here the $0$-neighbourhood of $t$ is simply $t$ itself. A $k$-collared
tile is a tile labelled with its $k$-neighbourhood. A $k$-collared
prototile is a translational congruence class of a $k$-collared tile.
Since $T$ has \flc\ there are only finitely many 
$k$-collared prototiles. The CW-complex $\Gamma^{k}$ is the complex
obtained from the disjoint union of all  $k$-collared prototiles upon
identifying any two faces $f_1,f_2$ of two $k$-collared prototiles
$p_1,p_2$ if one finds in $T$ two tiles $t_1,t_2$ whose 
$k$-neighbourhood corresponds to the two $k$-collared prototiles $p_1,p_2$
such that the faces of $t_1$ and $t_2$ which correspond to $f_1,f_2$
agree in $T$. There are surjective maps 
$\alpha_{k}:\Gamma^{k+1}\to\Gamma^{k}$
associating to a point in a $k\!+\!1$-collared prototile the same
point in the $k$-collared prototile obtained by reducing the labelling
of the $k\!+\!1$-neighbourhood to the $k$-neighbourhood. Then one has
the following result \cite{Gaehler,Sadun03}: $\Omega_T$ is the inverse limit of
the chain 
$$\cdots\stackrel{\alpha_{k+1}}{\longrightarrow}
\Gamma^{k+1}\stackrel{\alpha_{k}}{\longrightarrow} 
\Gamma^{k}\cdots \stackrel{\alpha_{0}}{\longrightarrow}
\Gamma^{0}.$$ 
Moreover we have maps $\pi_k:T\to\Gamma^k$ associating to a point $x$
in a tile of $T$ the point at the same position in the $k$-collared
prototile to which it corresponds.
  
As a consequence $\check H(\Om_\Pp,A)$ is the direct limit of the chain
$$H(\Gamma^{0},A) \stackrel{\alpha_{0}^*}{\longrightarrow}
\cdots H(\Gamma^{k},A) \stackrel{\alpha_{k}^*}{\longrightarrow}
H(\Gamma^{k+1},A)\cdots .$$ Here $H(\Gamma^{k},A)$ is the
cohomology of CW-complexes. In particular there are homomorphisms
$I_k: H(\Gamma^{k},A) \to \check H(\Om_\Pp,A)$ satisfying
$I_k=I_{k+1}\circ\alpha_k^*$.

\subsubsection{\eP\ cohomology}\label{sect-SK}
The \sP\ forms on $\rn$ form  
a sub-complex of the de Rham complex for $\rn$. The cohomology of this complex
is called the
strongly $P$-equivariant cohomology (of $\rn$) and denoted by 
$H_{s\dash \Pp}(\rn,\R)$. 
It has been shown in \cite{KellendonkPutnam06} that $H_{s\dash \Pp}(\rn,\R)$ is
isomorphic to $\check H(\Om_\Pp,\R)$. This isomorphism has been revisited
in \cite{Sadun} where the description of $\Omega_P$ as
inverse limit was used. We provide an explicit description of
this isomorphism in degree one.

Let $T$ be a polyhedral tiling of \flc\ and set $P=T^0$.
Recall the quotient map
$\pi_k:T\to\Gamma^k$. Let $e$ be a $1$-cell in $\Gamma^k$. Then
$\pi^{-1}_k(e)$ is a union of (oriented) edges of $T$. Suppose that
$d\varphi$ is \eP\ with range $r$ and that $r$ is small compared with $k$
(i.e. the tubular $r$-neighbourhood of any tile in $T$ is smaller than
its $k$-neighbourhood). Then $\int_{\tilde e}d\varphi$ is independent
of the choice of $\tilde e\in \pi^{-1}_k(e)$, since two such preimages
agree on their $k$-neighbourhood. Denoting by $Z^1_{r\!-\!P}(\rn,\rn)$ the
closed \eP\ $1$-forms of range $r$ we may thus define
$J_r^k:Z^1_{r\!-\!P}(\rn,\rn)\to Z^1(\Gamma^k,\rn)$ by 
\begin{equation}\label{eq-iso1}
J_r^k(d\varphi)(e) = \int_{\tilde e}d\varphi= 
\varphi(t(\tilde e))-\varphi(s(\tilde e)),
\end{equation}
where $s(\tilde e)$ and $t(\tilde e)$ are the vertices at which
$\tilde e$ starts and terminates.  
It is clear that $J_r^k(d\varphi)$ is closed, i.e.\ vanishes on
cycles. $J_r^k(d\varphi)$ is exact if and only if also 
$\varphi(s(\tilde e))$ is independent of the choice of the preimage
for $e$ which means that $\varphi$ is \sP\ (with a range with is
large compared $k$ but finite). It follows, that
$\Theta:H^1_{s\!-\!P}(\rn,\rn)\to \check H^1(\Omega_P,\rn)$:
$$ \Theta[d\varphi] = I_k[J_r^k(d\varphi)]$$
is well defined. It doesn't depend
on the precise choice for $r$ and $k$ as long as $d\varphi$ is \eP\ with
range $r$ and $k$ large compared to $r$.
Using this flexibility one easily sees that $\Theta$ is additiv and
hence defines a group homomorphism. 

We show that $\Theta$ is surjective: Let $f\in Z^1(\Gamma^k,\rn)$
choose $x_0\in T^0$ and define $\psi :T^0\to\rn$ by
\begin{equation}\label{eq-iso2} 
\psi(t(\tilde e))-\psi(s(\tilde e)) =
f(\pi_k(\tilde e)),\quad 
\psi(x_0)=0,
\end{equation}
for all (oriented) edges $\tilde e$ of $T$.  Let $h\in\rn$ be the
vector corresponding to the edge $\tilde e$. Then 
$\Delta_h\psi(x)  = \Delta_h\psi(y)$ for all $x,y\in
T^0\cap (T^0-h)$ which have the same $k+1$-neighbourhood. It follows
that $\Delta_h\psi$ is \eP\ with a range large compared to
$k+1$.
By Lemma~\ref{lem-extend} Part~\ref{item-3} there exists a smooth
function $\varphi$ which coincides with $\psi$ on $P$ and has \sP\
differential. It follows that for large enough $r$ and $l$,
$J_r^l(d\varphi) = f\circ \alpha_{l-1}\cdots\alpha_k$ and hence
$\Theta([d\varphi]) = I_k([f])$. 
Hence $\Theta$ is surjective. To see that $\Theta$ is injective we
note that by (\ref{eq-iso1}) $J_r^k(d\varphi)=0$ whenever $\varphi$ is
constant on the vertices of the tiling.
From Lemma~\ref{lem-trivial}
it follows then that $\varphi$ is \sP, hence $[d\varphi]=0$.
\bigskip

Also the \wP\ forms on $\rn$ form  
a sub-complex of the de Rham complex for $\rn$ and hence a cohomology
group which is called the
\wP\ cohomology of $\rn$. We will not consider \wP\ cohomology 
in this article, but one of our conclusions will be that
the mixed group
$$Z^1_{s\dash P}(\rn,\rn) / B^1_{w\dash P}(\rn,\rn) \cap Z^1_{s\dash
  P}(\rn,\rn)$$ 
is worth further investigation.

\section{The first cohomology and deformations}
\subsection{The approach of Williams, Sadun and Clark,
  Sadun}\label{sect-CS} 
Let $T$ be a polyhedral tiling of \flc\ with its system of 
Anderson-Putnam-G\"ahler complexes $(\cp^k)_k$. 

A shape function is a $\rn$-valued $1$-cocycle on $\cp^k$, i.e.\ a function
$f:\cp^k \to\rn$ such that $\delta f(w) := f(\delta w) = 0$ 
for any $2$-cell $w$ in $\cp^k$ ($\delta w$ is the chain
corresponding to the boundary of $w$).
With a little luck a shape function
defines a new tiling (up to translation) which is obtained from $T$ by
changing the shape of the prototiles. A $1$-cell $e$ of 
$\Gamma^k$ is an (oriented) egde of a $k$-collared prototile.
$f(e)$ is a vector in $\rn$
and hence defines an oriented straight line segment up to translation
which is going to be an edge of a new prototile.
Replacing each edge $e$ of a $k$-collared prototile 
by a line segment parallel to $f(e)$, the condition $\delta f=0$
guaranties that these line segments can be matched together at their
boundary points
with the same combinatorics as in the prototile.
One thus obtaines a new $1$-skeleton. Apriori there is no reason that
one can complete 
the new $1$-skeleton to a tile by adding in
higher dimensional faces in such a way that the combinatorial structure
is that of the prototile. But if this is the case, we can construct
from $T$, tile by tile, a whole new tiling up to an overall translation, and
hence a new tiling space which we denote by $(\Omega_T)_f$. 
It is somewhat clear that all this works if the new edges differ
only slightly from the old edges. It then makes sense to speak of a
deformation of $T$.

A question which we will investigate in more detail than has been done in
\cite{ClarkSadun06} is whether the above
process is invertible, i.e.\ whether one can obtain $T$ back from the new
tiling by a similar procedure. 
Certainly, if the new tiling would be periodic whereas the old
one was not, the process could not be inverted. One may argue that 
invertibility should be no problem if the new edges differ
little from the old edges. This is, however, only obvious if
$k$ is fixed, since the larger $k$ the more shapes are redefined by
the shape function.

Let us call shape functions which allow for all this admissible.
\begin{thm}[\cite{SadunWilliams03,ClarkSadun06}]\label{thm-sw}
Let $f: \cp^0 \to\rn$ be a cocycle which is admissible. Then $\Omega_T$
is homeomorphic to $(\Omega_T)_f$. 
\end{thm}
\begin{thm}[\cite{ClarkSadun06}]\label{thm-cs1}
Let $f,g: \cp^0 \to\rn$ be two cocycles which are admissible. 
If $I_0(f)=I_0(g)$ then $(\Omega_T)_f$
is mutually locally derivable with $(\Omega_T)_g$. 
\end{thm}
Here the relation of being mutually locally derivable is extended to
hulls of tilings by saying that two such hulls are  mutually locally
derivable if there exists a topological conjugacy between them which
maps a tiling to a mutually locally derivable tiling with a uniform
derivability range.

Clark and Sadun introduce the notion of asymptotically negligible
elements in $H^1(\Omega_T,\rn)$. Let us suppose that $r$ is larger
than the maximal diameter of a prototile. Then an ordered pair of 
vertices $(x_1,x_2)$ of $T$ is called a recurrence of size $r$ if
$B_{r}[T-{x_1}] =  B_{r}[T-{x_2}]$ but $B_{r+\delta}[T-{x_1}]
\neq B_{r+\delta}[T-{x_2}]$ for $\delta>0$. 
If $k$ is small
enough so that the $k+1$-neighbourhood of any tile is
contained in the $r$-tubular neighbourhood of the tile, then a path
between $x_1$ and $x_2$ along edges defines a loop in the $1$-skeleton
of $\Gamma^{k}$, different paths leading to homologous loops. 
An element $\eta\in H^1(\cp^k,\rn)$ can be evaluated on (the cycle
defined by) a loop in $\cp^k$.
It is called asymptotically negligible if 
for all $\epsilon$ there exists an $r_\epsilon$ such that, when $\eta$
is evaluated on a loop in the $1$-skeleton
of $\Gamma^{k}$ which is defined by a recurrence of size larger than
$r_\epsilon$, then the result is smaller than $\epsilon$ in norm. 
The images under the $I_k$ of asymptotically negligible elements in
$H^1(\cp^k,\rn)$ are by definition the asymptotic elements of
$H^1(\Omega_T,\rn)$. 
\begin{thm}[\cite{ClarkSadun06}]\label{thm-cs2}
Let $f,g: \cp^0 \to\rn$ be two cocycles which are admissible. 
If $I_0(f)-I_0(g)$ is asymptotically negligible then $((\Omega_T)_f,\rn)$
is topologically conjugate to $((\Omega_T)_g,\rn)$. 
\end{thm}
It is rather clear that the three theorems extend to shape functions
on $\Gamma^k$. 
Note also that, if $f:\cp^k \to\rn$ is a shape function then
$f\circ\alpha_k:\cp^{k+1} \to\rn$ is also a shape function and it leads to
the same deformation.

\subsection{Negligible elements versus \wP\  $1$-forms}

The aim is to reinterprete negligible elements as  \wP\  $1$-forms and
to interprete Theorem~\ref{thm-cs2} in terms of a mixed cohomology
group.
\begin{thm}\label{thm-neg}
An element $\eta\in \check H^1(\Omega_P,\rn)$ is asymptotically
negligible if and only if $\Theta^{-1}(\eta)$ belongs to 
$B^1_{w\!-\!P}(\rn,\rn)\cap Z^1_{s\!-\!P}(\rn,\rn)/B^1_{s\!-\!P}(\rn,\rn)$.
\end{thm}
\begin{proof}
Suppose $d\varphi\in B^1_{w\!-\!P}(\rn,\rn)\cap
Z^1_{s\!-\!P}(\rn,\rn)$. Hence for all $\epsilon>0$ exists
$R_\epsilon$ such that $\varphi=\varphi_\epsilon+\psi_\epsilon$ where
$\varphi_\epsilon$ is \eP\ with range $R_\epsilon$ and
$\|\psi_\epsilon\|_\infty <\epsilon$. Also $d\varphi$ is \eP\ with some finite
range which we may suppose to be smaller than $R_\epsilon$.
Let $k$ be large compared to $R_\epsilon$ and $r_\epsilon$ be large
compared to $k+1$. Finally let $(x_1,x_2)$ be a recurrence of size
larger or equal than $r_\epsilon$. Denote by $[x_1,x_2]$ the cycle it
defines in $\Gamma^k$.
Then we have
$$J_{R_\epsilon}^k(d\varphi)([x_1,x_2]) = 
\varphi_\epsilon(x_2)-\varphi_\epsilon(x_1) + 
\psi_\epsilon(x_2)-\psi_\epsilon(x_1).$$ The first difference drops out 
since $[x_1,x_2]$ is a loop in $\Gamma^k$ and so
$|J_{R_\epsilon}^k(d\varphi)([x_1,x_2])|<2\epsilon$. This shows that
$\Theta([d\varphi])$ is asymptotically negligible.

For the converse suppose that $f\in Z^1(\Gamma^k,\rn)$ such that $[f]$ 
is asymptotically negligible. As explained in
Section~\ref{sect-SK}, $\Theta^{-1}([f])$ is represented by the differential
of a smooth \sP\ function $\varphi:\rn\to\rn$
which extends the function
$\psi :T^0\to\rn$ defined as in (\ref{eq-iso2}).
Let $\epsilon>0$ and $r_\epsilon$ such that $f$ applied to
a loop defined by a recurrence  $(x_1,x_2)$ of size larger than or
equal to $r_\epsilon$ is smaller in norm than $\epsilon$. 
By the \flc\ of $T$ 
there exist finitely many points $p_1,\cdots,p_f\in T^0$ such that
the sets $T^0_i:=\{y\in T^0|B_{r_\epsilon}[T-y]
=  B_{r_\epsilon}[T-p_i]\}$ are pairwise disjoint, locally derivable
from $T^0$, and their union exhausts $T^0$. Furthermore,
the elements of $\bigcup_i T^0_i\times T^0_i$ are precisely the
recurrences of size larger than or equal to $r_\epsilon$.
Hence for all $x_1,x_2\in T_i^0$ we have
$|\varphi(x_1)-\varphi(x_2)|=|f([x_1,x_2])|<\epsilon$. 
From Prop.~\ref{prop-set} follows therefore
that $\varphi$ is \wP.
\end{proof}

\section{Deformations and $P$-equivariant functions}
The aim of this section is to
present a theory of deformations which is formulated
entirely in the
framework of \eP\ functions. We will recover the results of
Theorems~\ref{thm-sw},\ref{thm-cs1},\ref{thm-cs2}
(see Theorems~\ref{thm-closed},\ref{thm-exact} and Cor.~\ref{cor-back2})
but are able to go further as we obtain 
refined information about the homeomorphisms involved.
As a consequence we are able to decide
whether two pattern spaces are homeomorphic (with homeomorphism
preserving the abstract transversal) by means of pattern
equivariant functions, and to a large extend this turns out to be the
case for topological conjugacy as well. 
Finally we investigate the question of
invertibility of deformations.

Our first aim is to interprete the elements entering in the first
$P$-equivariant cohomology with values in $\rn$.
Since $\rn$ is contractible, any closed $1$-form on $\rn$, 
$P$-equivariant or not, can be written as
$d\varphi$ for a smooth function $\varphi:\rn\to\rn$. 
The first $P$-equivariant cohomology is therefore related to
functions, but the restriction that their differential 
is $P$-equivariant does not imply that they are themself \eP\
which is what makes the cohomology non-trival.
For instance $\varphi=\id$ is never $P$-equivariant,
but $d\,\id$ is.

One may view $\Omega_P$ as the compactification of $\rn$ defined by
$P\in\rn$. It is then natural to ask how does a function 
$\varphi:\rn\to\rn$ affect the compactification, in
other words, how does $\Omega_{\varphi(P)}$ compare to $\Omega_P$? 
Furthermore, one may ask how the 
pattern dynamical systems compare.

We will often need to make the following assumption
about $\varphi$.
\begin{hypo}\label{hypo} 
For all $r>0$ exists $r'>0$ such that for all $x\in \rn$: 
$$B(\varphi(x),r)\cap\varphi(P)\;\; \subset \;\;
\varphi(B(x,r')\cap P)\, .$$ 
\end{hypo}
This assumption is for instance verified if 
\begin{equation}\label{eq-H}
\forall r>0 \exists
r'>0\forall x\in \rn: \quad \varphi^{-1}(B(\varphi(x),r))\subset
B(x,r').
\end{equation}
Note that if $\varphi$ satisfies (\ref{eq-H})
and $\eta$ is bounded then also $\varphi-\eta$ satisfies (\ref{eq-H}).
We say that  
$\varphi:\rn\to\rn$ is a bi-Lipschitz map if
there exists $\lambda>1$ (the Lipschitz constant) such that
$$\lambda^{-1}|x-y| \leq |\varphi(x)-\varphi(y)|\leq \lambda |x-y|.$$
A bi-Lipschitz map is a homeomorphism which satisfies  (\ref{eq-H})
with $r'=\lambda r$. In fact, 
it is injective and continuous, and hence has open image by Brouwer's theorem.
Moreover its image is also closed since a pre-image
of a Cauchy sequence is a Cauchy sequence, and hence it is surjective.
A further advantage of a bi-Lipschitz
map is that it maps Delone sets to Delone sets.
We are particularily interested in the case 
$\varphi = \id +\eta$ where $\eta$ is differentiable and
$\|d\eta\|_\infty< 1$ 
(or $\varphi = g +\eta$ with $g\in GL(n,\R)$ and
$\|d\eta\|_\infty<(\|g^{-1}\|_\infty)^{-1}$). 
Such a map $\varphi$ is a bi-Lipschitz map 
with $\lambda^{-1}=1-\|d\eta\|_\infty$.

\begin{lemma}\label{lem1}
Let $P$ be a uniformly discrete set of \flc.
Let $\varphi:\rn\to\rn$ be a function,
$g\in GL(n,\R)$ and $\varphi-g$ be \eP\ with range $R$.
Then we have for all $r>0$ and $x,y\in\rn$:
\begin{eqnarray*}
&&B_{R+r}[P-x] = B_{R+r}[P-y] \Longrightarrow \\
&&\qquad \forall h\in
B(0,r):\, \varphi(x+h) - g(x) = \varphi(y+h) - g(y)\, .
\end{eqnarray*}
In particular
$$ B_{R+r}[P-x] = B_{R+r}[P-y] \Longrightarrow
B_r(\varphi(P) - g(x)) =
B_r(\varphi(P) - g(y)) $$
and $\varphi(P)$ is locally derivable from $g(P)$.
\end{lemma}
\begin{proof} By assumption
$B_{R+r}[P-x] = B_{R+r}[P-y]$ implies $\forall h\in
B(0,r):\, \varphi(x+h) - g(x+h) = \varphi(y+h) - g(y+h)$. The
statement follows therefore from the linearity of $g$. The last statement
is trivial for $g=\id$. Now if $\varphi-g$ is
\sP\ then $\varphi\circ g^{-1}-\id$ is
strongly $g(P)$-equivariant. This implies the last statement for general $g$.
\end{proof}

Note that  $d(\varphi-g)$ is \sP\ whenever $d\varphi$ is \sP.
\begin{lemma}\label{lem2} 
Let $P$ be a uniformly discrete set of \flc\
and $\varphi:\rn\to\rn$ be a differentiable function.
If $d\varphi$ is $P$-equivariant with range $R$
then for all $x,y\in\rn$:
\begin{eqnarray*}
&&B_{R+r}[P-x] = B_{R+r}[P-y] \Longrightarrow \\
&&\qquad \forall h\in
B(0,r):\, \varphi(x+h) - \varphi(x) = \varphi(y+h) - \varphi(y)\, .
\end{eqnarray*}
If moreover $\varphi$ satisfies Hypothesis~\ref{hypo} then for all $r>0$
exists $r'$ such that
\begin{equation}\label{comder}
B_{r'}[P-x] = B_{r'}[P-y] \Longrightarrow 
B_r(\varphi(P) - \varphi(x)) = B_r(\varphi(P) -  \varphi(y))\, .
\end{equation}
\end{lemma}
\begin{proof}
Let $B_{R+r}[P-x] = B_{R+r}[P-y]$ and $h\in B_r$. 
As $d\varphi$ is $P$-equivariant
with range $R$ we have $d\varphi(x+h)=d\varphi(y+h)$. In particular 
$\int_x^{x+h} d\varphi=\int_y^{y+h} d\varphi$ (we can integrate along
a straight line). As $\varphi(x+h)-\varphi(x)=\int_x^{x+h} d\varphi$ 
the first statement holds true.

Suppose  $B_{R+r'}[P-x] = B_{R+r'}[P-y]$. According to the 
first statement this implies
\begin{equation}\label{eq-H1}
\varphi(B(x,r')\cap P) -\varphi(x) =\varphi(B(y,r')\cap
P)-\varphi(y).
\end{equation} 
Now let $r$ be given and choose $r'$ according to
Hypothesis~\ref{hypo}.
Then, for any $z\in\rn$, 
$(\varphi(B(z,r')\cap P) -\varphi(z)) \cap
B(0,r) = (\varphi(B(z,r')\cap P) \cap B(\varphi(z),r)) -\varphi(z)
= (\varphi(P) \cap B(\varphi(z),r)) -\varphi(z)$. The latter is by
definition $B_r[\varphi(P)-\varphi(z)]$ and so (\ref{eq-H1})
implies $B_r[\varphi(P)-\varphi(x)]=B_r[\varphi(P)-\varphi(y)]$.
\end{proof}

\begin{cor}
Let $P$ be a uniformly discrete set of \flc.
Let $\varphi:\rn\to\rn$ be a differentiable function which fullfills
 Hypothesis~\ref{hypo}. If $d\varphi$ is strongly $P$-equivariant then  
$\varphi(P)$ has finite local complexity.
\end{cor}
\begin{proof}
Finite local complexity of $P$ means that there are
 for any given $r'$ only finitely many different  $B(0,r')\cap(P-x)$,
 $x\in P$. Thus we conclude from Lemma~\ref{lem2} that there are also only  
finitely many different $B(0,r)\cap
(\varphi(P)-z)$, $z\in\varphi(P)$.
\end{proof}
There is no reason for $\varphi(P)$ to be of \flc\ if we require merely
weak $P$-equivariance of $d(\varphi-g)$ or $d\varphi$.

How can we interprete (\ref{comder})? Clearly it allows to construct
the $r$-patch of $\varphi(P)$ around $\varphi(x)$ from the $r'$-patch
of $P$ around $x$, but this involves the passage from
$\varphi(x)$ to $x$ which is apriori not local (i.e.\ not determined
by a patch of $P$). 
Nevertheless, if we
forget about the position of the patches in $\rn$ and consider their
translational congruence classes we obtain from (\ref{comder}) a local
map: to each translational congruence class of an $r$-patch of
$P$ is associated a unique translational congruence class of an $r'$-patch of
$\varphi(P)$.
In one dimension this gives a rewriting rule. In higher dimensions one
obtains a local derivation in the sense of \cite{Kellendonk97}. If
$\varphi=\id+\eta$ with small $\|d\eta\|_\infty$ and $P$ is the set
of vertices of a polyhedral tiling the local map defined by (\ref{comder})
can be interpreted as a deformation 
of the tiling induced by changing
the length of the edges of its tiles \cite{SadunWilliams03}, see
Sect.~\ref{sect-CS}. In view of this and the results in Sect.~\ref{sect-SK}
we define:
\begin{defn} Let $P$ be a Delone set in $\rn$.
A Delone set $P'\subset\rn$ is a deformation of $P$ if there exists a
smooth function $\varphi:\rn\to\rn$ satisfying Hypothesis~\ref{hypo}
and having a \sP\ differential such that $P'=\varphi(P)$. 
\end{defn} 
We have a lot of freedom to choose $\varphi$. In fact,
since the joint kernel of the maps $J_r^k$ from Sect.~\ref{sect-SK} is
$$N^1_{s\dash P}(\rn,\rn)=\{d\eta\in B^1_{s\dash P}(\rn,\rn)|\eta
\mbox{ is constant on } P\}$$ 
the elements of $Z^1_{s\dash P}(\rn,\rn)/N^1_{s\dash P}(\rn,\rn)$
parametrise the deformations of $P$ (or of a tiling $T$ with $P=T^0$).

We end this section with the following important definition, which is
inspired from the above lemmata.
\begin{defn} 
Let $P$ be a uniformly discrete set of \flc.
\begin{enumerate}
\item Let $\varphi:\rn\to\rn$ be
a $C^0$-function such that $\varphi-g$ is \wP. We define 
$\rphi_g : \orb(P)\to\Omega_{\varphi(P)}$,
$$ \rphi_g(P-x) := \varphi(P) - g(x). $$
\item
Let $\varphi:\rn\to\rn$ be
a differentiable function such that $d\varphi$ is \sP.
We define
$\tphi: \orb(P)\to\Omega_{\varphi(P)}$,
$$ \tphi(P-x) := \varphi(P) - (\varphi(x)-\varphi(0)).$$ 
\end{enumerate}
\end{defn}
That these maps are well-defined is the issue of the following lemma.
\begin{lemma}: We assume the notation and conditions of the last definition.
In the first case
$P-x=P$ implies  $\varphi(P) - g(x)=\varphi(P)$.
In the second case
 $P-x=P$ implies  $\varphi(P) - \varphi(x)=\varphi(P)-\varphi(0)$.
\end{lemma}
\begin{proof}
Suppose that $\varphi-g$ is \sP. 
By Lemma~\ref{lem1} $\forall h\in P$:
$\varphi(-x+h)-g(-x) = \varphi(h)$ which implies
$\varphi(P)=\varphi(P)-g(x)$. 
If $\varphi-g$ is merely \wP\ we can find
for all $\epsilon>0$ an approximation
$\varphi_\epsilon$ of $\varphi$ such that $\varphi_\epsilon-g$ is
\sP\  and $\|\varphi_\epsilon-\varphi\|_\infty<\epsilon$.
The above then implies
that  $d_h(\varphi(P) - g(x),\varphi(P))<2\epsilon$. 

Now suppose that $d\varphi$ is \sP.
By Lemma~\ref{lem2} $\forall h\in P$:
$\varphi(x+h)-\varphi(x) = \varphi(h)-\varphi(0)$. This implies
$\varphi(P+x) = \varphi(P) +\varphi(x)-\varphi(0)$ and hence the statement.
\end{proof}

In the following we shall see that $\rphi_g$ and $\tphi$ can
be extended to maps between $\Omega_P$ and $\Omega_{\varphi(P)}$, the
second only in the \sP\ case. This
is the case once we can show that they are uniformly continuous, since a
uniformly continuous map between (not neccessarily complete) metric
spaces has a unique continuous extension to the completions.

\subsection{Closed strongly $P$-equivariant $1$-forms and
deformations}\label{sect-4.2}
We investigate more closely $P$-equivariant $1$-forms, that is \eP\ 
differentials, and how they affect the tiling spaces and dynamical systems.
We start by the observation that a \sP\ $1$-form gives rise to a
continuous map between hulls of Delone sets preserving the canonical
transversal, an observation which is slightly stronger than 
Theorem~\ref{thm-sw}. 
\begin{thm}\label{thm-closed}
Let $P$ be a Delone set of \flc.
Let $\varphi:\rn\to\rn$ be a $C^1$-map which fullfills
Hypothesis~\ref{hypo}
and whose differential $d\varphi$ is \sP. Then  
$\tphi$ is uniformly continuous map.
It thus extends  to a continuous map 
$\tphi:\Omega_P\to\Omega_{\varphi(P)}$.
If $\varphi(0)=0$ then $\tphi(\Xi_P) = \Xi_{\varphi(P)}$.
If $\varphi$ is surjective then
$\tphi$ is surjective.
\end{thm}
\begin{proof} 
Choose $r>0$ and then $r'\geq r$ according to Lemma~\ref{lem2}.
Let $D(P-x,P-y)<\frac{1}{r'}$. Since $P$ has \flc\ we can
  find $x',y'\in B(0,\frac{1}{2r'})$ such that 
$B_{r'}[P-x-x'] = B_{r'}[P-y-y']$. By Lemma~\ref{lem2} we then have
$B_{r}[\varphi(P)-\varphi(x+x')] = B_{r}[\varphi(P)-\varphi(y+y')]$ and so 
$D(\varphi(P)-\varphi(x+x'),\varphi(P)-\varphi(y+y'))<\frac{1}{r}$.
Since $P$ has \flc\ $d\varphi$ is bounded and
$|\varphi(x+x')-\varphi(x)|\leq \|d\varphi\|_\infty|x'|$.  
If $r$ is sufficiently large
then $D(\varphi(P)-\varphi(x+x'),\varphi(P)-\varphi (x))\leq
\|d\varphi\|_\infty|x'|  $. 
Hence the triangle inequality gives
$$
D(\varphi(P)-\varphi(x),\varphi(P)-\varphi(y)) \leq
\frac{1}{r}+\frac{\|d\varphi\|_\infty}{2r'}+\frac{\|d\varphi\|_\infty}{2r'}
\, .   $$
This tending to $0$ for $r$ tending to infinity we conclude that
$\tphi$ is uniformly continuous. 

That $\tphi$ preserves the canonical transversal provided $\varphi(0)=0$
follows immediately from the construction of $\tphi$.
Since $\tphi(\Xi_P)$ is closed and contains
$\{\varphi(P)-\varphi(p)|p\in P\}$ we must have $\tphi(\Xi_P) =
\Xi_{\varphi(P)}$. 

Finally, if $\varphi$ is surjective then $\orb(P')=\tphi(\orb(P))$
which is dense and hence implies surjectivity of $\tphi$.
\end{proof}
Note that since $d\varphi$ is \sP\ $\varphi(P)$
has \flc, but it need not be a Delone set. 
\bigskip

We are now interested in the converse direction: given a continuous
map between hulls of Delone sets, can we associate a \sP\ $1$-form to it?

Let $P,P'$ be aperiodic Delone sets of \flc\ and
$\Phi:\Omega_P\to\Omega_{P'}$ a continuous map satisfying $\Phi(P)=P'$. 
Since the orbits under the translation action are the path-connected
components $\Phi$ must preserve orbits. 
Hence for any $x$ exists a
$x'$ such that $\Phi(P-x)=P'-x'$. Since $P'$ is aperiodic
$x'$ is uniquely determined by $x$. Thus there
exists a family of functions $\varphi_{P-y}:\rn\to\rn$, $y\in\rn$ such that 
\begin{equation}
\label{eq-var}
\Phi(P-y-x) = P'-y-\varphi_{P-y}(x),
\quad \varphi_{P-y}(0)=0.
\end{equation}
\begin{lemma}\label{lem-closed-converse}
For all $y\in\rn$, $\varphi_{P-y}:\rn\to\rn$ is uniformly continuous.
Furthermore, for all $x\in\rn$, 
$\orb(P)\ni \omega \mapsto\varphi_\omega(x)\in\rn$ is uniformly continuous.  
\end{lemma}
\begin{proof}
To show that $\varphi_P$ is uniformly continuous choose
$0<\epsilon<\frac{r_{min}(P')}{3}$,
where $r_{min}(P')$ is the minimal distance in $P'$. 
By the uniform continuity of $\Phi$ there
exists $\delta$ such that for all $x_0\in\rn$, $|v|<\delta$ implies
\begin{equation}\label{eq-leaf}
D(\Phi(P-x_0),\Phi(P-x_0-v))<\epsilon.
\end{equation}
By the remarks at the end of Sect.~\ref{sec-2.1} the pre-image under
$w\mapsto \Phi(P-x_0)-w$ of the
path-connected component of $\Phi(P-x_0)$ in $U_\epsilon(\Phi(P-x_0))$
is contained in $B(0,\epsilon)$. Therefore 
$\varphi_{P-x_0}(B(0,\delta))\subset B(0,\frac{r_{min}(P')}{3})$. 

Furthermore, (\ref{eq-leaf}) implies that either
$|\varphi_{P}(x_0+v)-\varphi_{P}(x_0)|\leq\epsilon$ or
$|\varphi_{P}(x_0+v)-\varphi_{P}(x_0)|\geq
r_{min}(P')-\epsilon$. But the second possibility cannot hold since
$|\varphi_{P}(x_0+v)-\varphi_{P}(x_0)|=|\varphi_{P-x_0}(v)|\leq
\frac{r_{min}(P')}{3}$. This proves uniform continuity of $\varphi_P$. That of
$\varphi_{P-y}$ follows in the same way.

To show that $\omega\mapsto \varphi_\omega(y)$ is uniformly continuous on the
orbit of $P$ we choose $\epsilon>0$, small compared with $|y|^{-1}$. 
Note that $B_r[\omega]=B_r[\omega']$ implies
$B_{r-|x|}[\omega-x]=B_{r-|x|}[\omega'-x]$ provided $r>|x|$. 
Using this property, the (uniform) continuity of $\Phi$
and boundedness of $\varphi_{\omega'}$ we can 
find $\delta>0$ such that  
$D(\omega,\omega')<\delta$ implies
for all $y'\in B(0,|y|)$: 
$D(\Phi(\omega-y'),\Phi(\omega'-y'))<\epsilon$ and
$D(\Phi(\omega)-\varphi_{\omega'}(y'),\Phi(\omega')-
\varphi_{\omega'}(y'))<\epsilon$ and therefore
$$
D(\Phi(\omega)-\varphi_{\omega'}(y'),\Phi(\omega)-
\varphi_{\omega}(y')) 
< 2\epsilon.
$$
This implies that $|\varphi_{\omega}(y')-\varphi_{\omega'}(y')|$
is either smaller than $2\epsilon$ or larger than
$r_{min}(P')-2\epsilon$. 
From the continuity
of $\varphi_\omega$ and the condition $\varphi_\omega(0)=0$ follows
therefore that $|\varphi_{\omega}(y)-\varphi_{\omega'}(y)|<2\epsilon$
(for small enough $\epsilon$).
\end{proof}
We now require in addition that $\Phi$ preserves the canonical transversal.
\begin{thm}\label{thm-closed-converse}
Let $P,P'$ be aperiodic Delone sets of \flc. Let
$\Phi:\Omega_P\to\Omega_{P'}$ be a continuous map satisfying
$\Phi(\Xi_P)\subset \Xi_{P'}$ and $\Phi(P)=P'$. Then 
there exists a smooth function $\varphi:\rn\to\rn$ with \sP\ differential
which coincides with $\varphi_P$ (as defined by (\ref{eq-var}))
on $P$ and satisfies $\varphi(P)\subset P'$.
\end{thm}
\begin{proof} We can sharpen the last part of the proof of
Lemma~\ref{lem-closed-converse} provided $w,w'\in \Xi_P$
and $y=P-P$.   
Let $\epsilon>0$, $M>0$ and $p,q\in P$.
By the reasoning of the proof of Lemma~\ref{lem-closed-converse} 
there exists $R>M$ such
that $B_R[P-p]=B_R[P-q]$ implies for all $y\in B_M[P-p]$:
$D(\Phi(P-p-y),\Phi(P-q-y))<\epsilon$ and
$|\varphi_{P-p}(y)- \varphi_{P-q}(y)|<\epsilon$.
Moreover, $\Phi(P-p-y)$ and $\Phi(P-q-y)$
lie in $\Xi_{P'}$ and hence
the finite local complexity of $P'$ 
allows us to replace $D$ by $D_0$ to obtain
$B_{\frac{1}{\epsilon}}[\Phi(P-p-y)]=B_{\frac{1}{\epsilon}}[\Phi(P-q-y)]$.
Hence
$|\varphi_{P-p}(y)-\varphi_{P-q}(y)|<\epsilon$ forces 
$|\varphi_{P-p}(y)-\varphi_{P-q}(y)|$ to vanish.
A simple calculation shows that
$\varphi_{P-p}(y)=\varphi_P(p+y)-\varphi_P(p)$. Hence for all $h\in P-P$,
$\Delta_h\varphi_P:P\cap (P-h)\to\rn$ is \sP. 
Applying
Lemma~\ref{lem-extend}.\ref{item-3} we obtain a smooth function $\varphi$ with
\sP\ differential which coincides with $\varphi_P$ on $P$. 

Finally, if
$p\in P$ then $\Phi(P-p)\in\Xi_{P'}$ and hence $0\in
P'-\varphi(p)$. This shows that $\varphi(P)\subset P'$. 
\end{proof}

We cannot expect that under the hypothesis of the theorem
$\varphi(P)$ is all of $P'$. For instance, 
if $P'$ is obtained from $P$ by adding
points in a locally derivable way then $\varphi=\id$ would be
the resulting map but $P$ a proper subset of $P'$.
\begin{cor}\label{cor-def} Let $P,P'$ be aperiodic Delone sets of \flc. Let
$\Phi:\Omega_P\to\Omega_{P'}$ be a
homeomorphism which maps $\Xi_P$ onto $\Xi_{P'}$ and $P$ to $P'$. 
Then the map
$\varphi$ constructed in the theorem satisfies Hypothesis~\ref{hypo} and
$P'=\varphi(P)$. 
\end{cor}
\begin{proof}
Since $\Phi$ is a homeomorphism also $\varphi_P$ is a homeomorphism
and its inverse $\varphi_P^{-1}$ is the map constructed in
(\ref{eq-var}) from $\Phi^{-1}:\Omega_{P'}\to \Omega_P$, we denote it
here by $\varphi_{P'}$. By Lemma~\ref{lem-closed-converse} it is
uniformly continuous and hence we can find at least one $r>0$ for
which $\exists r'>0,\forall y\in\rn$:  
$\varphi_{P'}(B(y,r))\subset B(\varphi_{P'}(y),r')$. 
Our aim is to
show that this implies that $\forall r>0,\exists r'>0,\forall y\in\rn$:  
$\varphi_{P'}(B(y,r))\subset B(\varphi_{P'}(y),r')$ and hence that
$\varphi_P$ satisfies (\ref{eq-H}).
To see this let $0\in A\subset B(0,r)$ be a finite set such that 
$B(0,\frac{3}{2}r)\subset \bigcup_{a\in A} B(a,r)$. Then the above
implies
$$\varphi_{P'}(B(y,\textstyle{\frac{3}{2}}r))\subset \bigcup_{a\in A}
B(\varphi_{P'}(y+a),r').$$
By Lemma~\ref{lem-closed-converse} 
the maps $y\mapsto \varphi_{P'}(y+a)-\varphi_{P'}(y)$ are weakly
$P'$-equivariant and hence bounded. If we denote a common bound for
their norm by $M$ then $\bigcup_{a\in A}
B(\varphi_{P'}(y+a),r')\subset B(\varphi_{P'}(y),r'+M)$.
This shows that $\varphi_P$ satisfies (\ref{eq-H}).
Since $\varphi_P$ and $\varphi$
are uniformly continuous ($\varphi$ as $d\varphi$ is bounded),
$\varphi_P-\varphi$ is bounded and hence $\varphi$ satisfies 
Hypothesis~\ref{hypo}.

We have $\Phi(\orb(P))=\orb(P')$ and
$\Phi^{-1}(\Xi_{P'}) = \Xi_{P}$. Now if $y\in P'$ then the first
identity implies that there exists a $x$ such that
$P'-\varphi_P(x)=P'-y$, and hence $\varphi_P(x)-y$ lies in the
periodicity lattice of $P'$. Using $\varphi(0)=0$ and
continuity we can arrange for $\varphi_P(x)=y$. Now the second
identity implies that $P-x\in\Xi_P$ and hence $x\in P$.
\end{proof}
Let us call two homeomorphisms from the hull of $P$ into the hull 
of another pattern $P'$ 
which preserve the canonical transversals and map $P$ to $P'$ 
equivalent if the maps constructed in (\ref{eq-var}) for $y=0$
coincide on $P$. Then we may say
that a small neighbourhood of the class of $d\id$ in
$Z^1_{s\dash P}(\rn,\rn)/ N^1_{s\dash P}(\rn,\rn)$ parametrises
equivalence classes of such homeomorphisms.

\begin{remark}
The map $\varphi$ constructed in Theorem~\ref{thm-closed-converse} 
may differ substantially from $\varphi_P$ away from $P$. 
A refinement of Lemma~\ref{lem-extend}.\ref{item-3} 
allows us to construct maps which are arbitrarily close, namely
for all $\epsilon>0$ 
there exists a smooth function $\varphi:\rn\to\rn$ with \sP\ differential
such that $\|\varphi_P-\varphi\|_\infty<\epsilon$ and $\varphi$ coincides
with $\varphi_P$ on $P$.
 
This can be seen as follows. Fix $\delta>0$. Uniform continuity of
$\omega\mapsto \varphi_\omega(x)$ is equivalent to
weak $P$-equivariance of 
$y\mapsto\varphi_{P-y}(x)=\Delta_x\varphi_P(y)$
($x\in\rn$ fixed).\footnote{We extend the
  definition of $\Delta_x\varphi_P(y)=\varphi_P(y+x)-\varphi_P(y)$ to
  $x,y\in\rn$.} 
Hence we can find a partition  $\{P_i\}_i$ of $P$ such that for all
$x\in B(0,3A)$ and $p,q\in P_i$ we have
$|\Delta_x\varphi_P(p)-\Delta_x\varphi_P(q)|<\delta$. 
Furthermore, taking if neccessary a refinement of that partition the
strong $P$-equivariance of $\Delta_h \varphi_P:P\cap(P-h)\to\rn$ allows us to
assume that for $h\in P-P, |h|\leq
A$, $\Delta_h\varphi_P$ takes the same value on all points of $P_i$. 
Choose $p_i\in P_i$ and set, for $p\in P_i$ and $x\in \mbox{\rm
int} V_p-p$,
\begin{equation}\label{eq-ext2}
\tilde\phi(x+p) = \psi(x+p_i)-\psi(p_i)+\psi(p) .
\end{equation}
Now we proceed exactly as in the proof of
Lemma~\ref{lem-extend}.\ref{item-3} 
thereby obtaining a smooth function $\varphi=\rho*\tilde\phi$ which has \sP\
differential. If we choose $\delta$ small enough then
$\|\varphi-\varphi_P\|_\infty$ is smaller than any given $\epsilon$.
The only drawback is that $\varphi-\varphi_P$ does not vanish on $P$.
But the restriction of $\varphi-\varphi_P$ to $P$ is \sP. We may therefore
use Lemma~\ref{lem-extend}.\ref{item-3} to add a smooth \sP\ function to
$\varphi$ such that the result coincides with $\varphi_P$ on $P$.
This correction is of order $\epsilon$ in the $\|\cdot\|_\infty$-norm
and so we are done.
\end{remark}

\subsection{Exact $P$-equivariant $1$-forms}

Exact $P$-equivariant $1$-forms are of the form $d f$ where $f$
is $P$-equivariant. This implies that $f$ is bounded. Such maps are
not good candidates to compare $P$ with $f(P)$ as the latter
could not even be uniformly discrete. It is more fruitful to consider 
exact one forms of the form $d(\varphi-\id)$ and then compare $P$ with
$\varphi(P)$.
\begin{thm}\label{thm-exact}
Let $P$ be a Delone set of \flc,
$\varphi:\rn\to\rn$ be a $C^0$-function and $g\in GL(n,\R)$. 
Suppose that $\varphi-g$ is \wP. Then  
$\rphi_g$ is a uniformly continuous map which satisfies
$\rphi_g(\omega-x)=\rphi_g(\omega)-g(x)$. 
In particular $\rphi_\fid$ 
extends to a
topological semi-conjugacy between $(\Omega_P,\rn)$ and 
$(\Omega_{\varphi(P)},\rn)$ mapping $P$ to $P'$.
If $\varphi$ is even \sP\ then $\varphi(P)$ is
a Delone set.
\end{thm}
\begin{proof} 
We consider first the case that $\varphi-g$ is \eP\ with range
  $R$. 
Clearly $g(P)$ is a Delone set. Since
$\varphi(P)=g(P)+(\varphi-g)(P)$ and $(\varphi-g)(P)$ is finite by
Lemma~\ref{lem-PE}, also $\varphi(P)$ is a Delone set. For the second
statement 
choose $r>0$ and let $r'=R+r$.
Let $D(P-x,P-y)<\frac{1}{r'}$. Since $P$ has \flc\ we can
  find $x',y'\in B(0,\frac{1}{2r'})$ such that 
$B_{r'}[P-x-x'] = B_{r'}[P-y-y']$. By Lemma~\ref{lem1} we then have
$B_{r}[\varphi(P)-g(x+x')] = B_{r}[\varphi(P)-g(y+y')]$ and so 
$D(\varphi(P)-g(x+x'),\varphi(P)-g(y+y'))<\frac{1}{r}$. 
If $r$ is sufficiently large
then $D(\varphi(P)-g(x+x'),\varphi(P)-g(x))\leq |g(x')| 
\leq \frac{\|g\|_\infty}{2r'}$. 
Hence the triangle inequality gives
$$
D(P-g(x),P-g(y)) \leq 
\frac{1}{r}+\frac{\|g\|_\infty}{R+r} \, .
$$
Since this tends to $0$ when $r$ tends to infinity we conclude that
$\rphi_g$ is uniformly continuous. 

Now suppose that $\varphi-g$ is merely \wP\ and hence
for any $\epsilon>0$ exists a \sP\ $C^0$-function
  $\varphi_\epsilon$ such that 
$\|\varphi-\varphi_\epsilon\|_\infty <\epsilon$.
Now 
$$ D(\rphi_g(P-x),{\rphi_\epsilon}_g(P-x))\leq 
d_H(\varphi(P),{\varphi_\epsilon}(P)) \leq \epsilon 
$$
and a further application of the triangle inequality implies that $\rphi_g$
is uniformly continuous.

The surjectivity of $g$ implies that the image of $\rphi_g$ contains
the orbit of $\varphi(P)$ and hence that  $\rphi_g$ is surjective.
The relation $\rphi_g(\omega-x)=\rphi_g(\omega)-g(x)$ is clear for
$\omega=P-y$ and follows for general $\omega$ from the
continuity of $\rphi_g$.
\end{proof}

We now investigate conditions which lead to statements 
in the converse direction: Assuming that $P'$ is a deformation of $P$
such that the dynamical systems are pointed topological conjugate, is
there a \wP\ map $\varphi$ such that $\varphi(P)=P'$?
Clark and Sadun obtained related results for substitution tilings 
\cite{ClarkSadun06}.

Consider two patterns $P$ and $P'$ of \flc\ and suppose that
$\Psi:\Omega_P\to\Omega_{P'}$ is a semi-conjugacy mapping
$P$ to $P'$. Since it commutes with the $\rn$ actions the maps 
constructed in Lemma~\ref{lem-closed-converse} are all equal to the identity.
Moreover, since 
$\Psi$ has to be the continuous extension of $ P-x \mapsto P'-x$ we
must have
$$\Psi=\rphi_\fid.$$
Of course, the semi-conjugacy is
not assumed to preserve the canonical transversals. 
If we assume in addition
that  $P'$ is the image of $P$ under a $C^1$-function $\varphi$ satisfying 
Hypothesis~\ref{hypo}
and having \sP\ differential we obtain from Theorem~\ref{thm-closed}
a (in general different) map $\tphi:\Omega_P\to\Omega_{P'}$ which
preserves the canonical transversal but does not commute with the
$\rn$ actions. Our main question is
therefore, under which additional conditions does
uniform continuity of $\rphi_\fid$ and $\tphi$
imply that $\varphi-\id$ is \wP?

A first observation to make is that boundedness of $\varphi-\id$ is
neccessary, as \wP\ functions are bounded. We show below that boundedness of
$\varphi-\id$ is sufficient, provided $P'$ is aperiodic and uniformly
continuous. 
\begin{lemma}\label{lem-aper}
Let $P\subset\rn$ be an aperiodic, uniformly discrete set.
For all $M>0$ there exists $r>0$ such that for all $h\in\rn, |h|\leq M$ holds:
if $B_r[P]=B_r[P-h]$ then $h=0$.
\end{lemma}
\begin{proof}
Let us suppose the contrary. Then there exist $M>0$,
a sequence $(r_n)_n$ tending to $+\infty$, and points 
$h_n\in \rn, |h_n|\leq M$ such that  $ B_{r_n}[P]=B_{r_n}[P-h_n]$ but
$h_n\neq 0$. 
Since $\{h\in\rn:|h|\leq M\}$ is compact we can assume, perhaps after
taking a sub-sequence, that $(h_n)_n$ converges, let's say to $h$.
This limit satisfies $h\neq 0$ since $B_{r_n}[P]=B_{r_n}[P-h_n]$ for
large enough $r_n$
implies that either $h_n = 0$ or $|h_n|\geq r_{min}(P)$.
It follows that $(P-h_n)_n$ converges to
$P-h$ and that $P-h = P$. This contradicts the
assumption that $P$ is aperiodic.
\end{proof}
\begin{thm}\label{thm-exact-converse}
Let $P'$ be an aperiodic and uniformly discrete 
deformation of a Delone set $P$, 
i.e.\ in particular $P'=\varphi(P)$
for a smooth function $\varphi:\rn\to\rn$ with \sP\ differential and
satisfying Hypothesis~\ref{hypo}.
Suppose that $P'$ is moreover pointed topologically 
semi-conjugate to $P$ and hence
$\rphi_\fid$ is uniformly continuous. If $\varphi-\id$ is bounded
then it is \wP.
\end{thm}
\begin{proof} Without loss of generality we may assume that $\varphi(0)=0$.
Then $\tphi:\Omega_{P}\to \Omega_{P'}$, which 
is uniformly continuous by Theorem~\ref{thm-closed}, is given by
$$\tphi(P-x) = \rphi_\fid(P-x) + \eta(x)$$
where $\eta = \varphi-\id$. By definition of the metric, uniform
continuity of $\tphi$ and $\rphi_\fid$ imply that for any
$\epsilon>0$ exists $\delta>0$ such that for all $x,y\in\rn$:
\begin{eqnarray*}
&& B_{\delta^{-1}}[P-x] =
B_{\delta^{-1}}[P-y] \\
\Longrightarrow && 
\left\{
\begin{array}{l}
B_{\frac{1}{\epsilon}}[\tphi(P-x)] \cong_{\epsilon}
B_{\frac{1}{\epsilon}}[\tphi(P-y)]  \\
B_{\frac{1}{\epsilon}}[\tphi(P-x)+\eta(x)] \cong_{\epsilon}
B_{\frac{1}{\epsilon}}[\tphi(P-y)+\eta(y)]
\end{array}
\right. 
\end{eqnarray*}
Suppose $\epsilon^{-1}> 2(\|\eta\|_\infty+1)$.  
Let $p\in P$. Then we have for all $x\in\rn$:
\begin{eqnarray*}\label{eq-5}
& & B_{\delta^{-1}}[P-x] =
B_{\delta^{-1}}[P-p] \Longrightarrow \exists x'\in\rn,|x'|< 2\epsilon:\\
& &
{}\quad B_{\frac{1}{2}\epsilon^{-1}}[\tphi(P-p)] =
B_{\frac{1}{2}\epsilon^{-1}}[\tphi(P-p)+\eta(x)-\eta(p)-x'].  
\nonumber
\end{eqnarray*}
Since $\eta$ is bounded we can apply Lemma~\ref{lem-aper} to
$\tphi(P-p)$ with 
$h=\eta(x)-\eta(p)-x'$ to conclude that $|\eta(x)-\eta(p)|<2\epsilon$.
By the \flc\ of $P$ there exists a finite subset 
$\{p_1,\cdots,p_k\}\subset P$ such that the sets
$\{p\in P|B_{\delta^{-1}}[P-p] = B_{\delta^{-1}}[P-p_i]\}$ partition $P$.
They are clearly locally derivable from $P$. Hence we can apply
Prop.~\ref{prop-set} to obtain the statement.
\end{proof}
The following example shows that the hypothesis of boundedness of
$\varphi-\id$ is not implied by the other conditions of the theorem.
The Penrose tilings
can be defined by means of a substitution.
This is a locally defined decomposition $\gamma$ of the tiles into smaller
tiles followed by a rescaling by factor of $\tau$, the golden
ratio. If $T$ is a Penrose tiling then $\gamma(T)$ and $T$ are mutually locally
derivable, since the substitution is recognisable. In
particular, $T$ and $\gamma(T)$ are pointed topologically conjugate
and so are their vertex sets $T^v$ and $\gamma(T)^v$.
Now a few Penrose tilings have the property that
they are invariant under the forth iterate of the substitution.
Hence for those, $\gamma^4(T)^v = \tau^{-4} T^v$, and thus $\gamma^4(T)^v$
is a deformation of $T^v$ given by the map $\varphi = \tau^{-4}\id$.
Clearly $\varphi-\id$ is unbounded. 

\subsection{Deforming back} $\varphi(P)$ being a deformation of $P$
does not imply that $P$ is a deformation of $\varphi(P)$, i.e.\
deformations are not always invertible (for instance if $\varphi(P)$
is periodic whereas $P$ was not).
Although this is not analysed in \cite{ClarkSadun06},
it is quite clear that, for fixed $k$ (cf.\ Section~\ref{sec-3.1}), 
there exists 
a neighbourhood of the identity shape function such that
shape functions in that neighbourhood are admissible and define
deformations which can be inverted. 
It is not clear apriori however, that one can obtain a lower
bound on the size of this neighbourhood which is uniform in $k$. 
We will in this section
effectively establish such a lower bound (Theorem~\ref{thm-back})
and thereby guarantee that
small deformations are invertible, i.e.\ can be deformed back.

Given a uniformly discrete set $P$ with minimal distance between
points $r_{min}$ we define
$A_P:[r_{min},\infty)\to \R^+$, 
$$A_P(r)=\inf\{|h-h'|:\,h, h'\in \bigcup_{x\in
  P}B_{r}[P-x], h\neq h'\}.$$
$A_P(r)$ is a monotonically decreasing
function which is bounded from above by $r_{min}$.
Moreover, if $P$ has \flc\ then $A_P(r)$ is strictly positive since then
$\bigcup_{x\in  P}B_{r}[P-x]$ is a finite set.

For the following lemma and its corollary we do not require any
$P$-equivariance property of $\varphi$ or $d\varphi$. 
In particular, $\varphi(P)$ will usually not have \flc.
\begin{lemma}\label{lem-back}
Let $P$ be a uniformly discrete set of \flc.
For all
  $r>0$ exists $1>\epsilon>0$ such that for all differentiable
  functions $\varphi:\rn\to\rn$ which satisfy
$\|d(\varphi-\id)\|_\infty<\epsilon$ holds:
$\forall x,y\in P$:
\begin{eqnarray*}
&& B_{r'}[\varphi(P)-\varphi(x)]=B_{r'}[\varphi(P)-\varphi(y)]\Longrightarrow\\
&&\qquad \left\{
\begin{array}{l}
B_{r}[P-x]=B_{r}[P-y]\quad\mbox{\rm and} \\
\forall h\in B_r[P-x]: \varphi(x+h)-\varphi(x)=\varphi(y+h)-\varphi(y)
\end{array}\right.
\end{eqnarray*}
where $r'=r(1-\epsilon)^{-1}$. 
\end{lemma}
\begin{proof}
The statement is trivial for $r<r_{min}$, the minimal distance in $P$,
so we assume $r\geq r_{min}$. Define
\begin{equation}\label{eps}\nonumber
\epsilon = \sup_{0<t < 1}\{2 t(1-t)^{-2 }r \leq A_P((1-t)^{-2 }r)\}.
\end{equation}
Suppose that $\|d(\varphi-\id)\|_\infty<\epsilon$.
Let $r'=r(1-\epsilon)^{-1}$ and 
$S_{r'}(x) :=\{h\in P-x | \varphi(x+h)-\varphi(x) \in
B_{r'}[\varphi(P)-\varphi(x)]\}$. 
Since $\varphi$ is bi-Lipschitz with constant $\lambda=
(1-\epsilon)^{-1}$ we have
\begin{equation}\label{incl}
B_{(1-\epsilon)r'}[P-x] \subset S_{r'}(x) \subset
B_{(1-\epsilon)^{-1}r'}[P-x].
\end{equation}
Now let $B_{r'}[\varphi(P)-\varphi(x)]=B_{r'}[\varphi(P)-\varphi(y)]$.
This implies
that there exists a bijection $\beta: S_{r'}(x)\to S_{r'}(y)$ such
that 
$$\varphi(y+\beta(h))-\varphi(y) = \varphi(x+h)-\varphi(x).$$ 
It follows that $ \beta(h)-h=
\eta(x+h)-\eta(x)) - \eta(y+\beta(h))-\eta(y)$ 
with $\eta=\varphi-\id$ and hence
$$ |\beta(h)-h|\leq (|h|+|\beta(h)|)\|d(\varphi-\id)\|_\infty .$$
By the right inclusion of (\ref{incl})
$$(|h|+|\beta(h)|)\|d(\varphi-\id)\|_\infty
< 2(1-\epsilon)^{-1}r'\epsilon = A_P((1-\epsilon)^{-1}r') .$$
Since $\beta(h),h\in \bigcup_{x\in P}B_{r}[P-x]$ 
the definition of $A_P$ implies that $h=\beta(h)$ and so
$S_{r'}(x)=S_{r'}(y)$.
From the left inclusion of (\ref{incl}) follows that
$B_{r}[P-x]=B_{r}[P-y]$
and, since $h=\beta(h)$,
 $\forall h\in B_r[P-x]: \varphi(x+h)-\varphi(x)=\varphi(y+h)-\varphi(y)$.
\end{proof}
\begin{cor}\label{cor-back}
Let $P$ be a Delone set of \flc. 
There exists $1>\epsilon>0$ such that for all differentiable
  functions $\varphi:\rn\to\rn$ which satisfy
$\|d(\varphi-\id)\|_\infty<\epsilon$ holds:
$\forall r>0$ $\exists r'>0$ $\forall x,y\in P$:
\begin{equation}
B_{r'}[\varphi(P)-\varphi(x)]=B_{r'}[\varphi(P)-\varphi(y)]\Longrightarrow
B_{r}[P-x]=B_{r}[P-y].
\end{equation}
\end{cor}
\begin{proof}
Choose $\delta>0$ and $r_0>0$ such that for all $r\geq r_0$ and
$x\in P$ holds
\begin{equation}\label{back1} 
B_{r+\delta}[P-x] \subset \bigcup_{h\in  B_{r}[P-x]}
B_{r}[P-x-h]+h.
\end{equation}
This is possible since $P$ is relatively dense.
Choose $\epsilon$ such that the statement of 
Lemma~\ref{lem-back} is satisfied with $r=r_0$ and some $r'=r_0'$.
Let, for $x\in P$,
$$ V(x):= \bigcup_{h\in  B_{r}[P-x]}
B_{r_0'}[\varphi(P)-\varphi(x+h)]+\varphi(x+h)-\varphi(x).$$
$ V(x)$ is a subset of $\varphi(P)-\varphi(x)$.
Then 
\begin{equation}\label{back2}
B(\varphi(x+h)-\varphi(x),r_0')\cap V(x) =
B_{r_0'}[\varphi(P)-\varphi(x+h)]+\varphi(x+h)-\varphi(x).
\end{equation}
Now let $V(x) = V(y)$, $x,y\in P$. This implies in particular   
$B_{r_0'}[\varphi(P)-\varphi(x)]=B_{r_0'}[\varphi(P)-\varphi(y)]$ 
from which we conclude from Lemma~\ref{lem-back} that
$B_{r_0}[P-x]=B_{r_0}[P-y]$ and $\forall h\in B_{r_0}[P-x]$:
$\varphi(x+h)-\varphi(x)=\varphi(y+h)-\varphi(y)$. The latter implies
with (\ref{back2}) that $\forall h\in B_{r_0}[P-x]$:
$B_{r_0'}[\varphi(P)-\varphi(x+h)] = B_{r_0'}[\varphi(P)-\varphi(y+h)]$
and hence by Lemma~\ref{lem-back} that
$\forall h\in B_{r_0}[P-x]$ holds $B_{r_0}[P-x-h]=B_{r_0}[P-y-h]$
and $\forall h'\in B_{r_0}[P-x-h]$: 
$\varphi(x+h+h')-\varphi(x+h)=\varphi(y+h+h')-\varphi(y+h)$.
In particular we find $\forall h''\in B_{r_0}[P-x-h]+h$: 
\begin{eqnarray*}
\varphi(x+h'')-\varphi(x)&=&\varphi(x+h+(h''-h))-
\varphi(x+h)+\varphi(x+h)-\varphi(x)\\
&=& \varphi(y+h+(h''-h))-
\varphi(y+h)+\varphi(y+h)-\varphi(y)\\
&=& \varphi(y+h'')-\varphi(y),
\end{eqnarray*}
as $ h''-h \in B_{r_0}[P-x-h]$.
Thus $V(x) = V(y)$, $x,y\in P$ implies 
$B_{r_1}[P-x]=B_{r_1}[P-y]$ and 
$\forall h\in B_{r_1}[P-x]$: 
$\varphi(x+h)-\varphi(x)=\varphi(y+h)-\varphi(y)$ where $r_1=r_0+\delta$.
Since $V(x)$ is finite there exists $r_1'$ such that $V(x)\subset
B(0,r_1')$. It follows that, for $x,y\in P$,
$B_{r_1'}[\varphi(P)-\varphi(x)]=B_{r_1'}[\varphi(P)-\varphi(y)]$
implies  $V(x) = V(y)$. We have thus proven that the statement of
 Lemma~\ref{lem-back} holds with the same $\epsilon$ but $r=r_1$ and $r'=r_1'$.
Iterating this argument yields the
statement of the corollary.
\end{proof}
Note that in the statement of the corollary the points $x$ and $y$
were taken 
from $P$ and so the corollary does not imply that $d\varphi^{-1}$ is
\sP.\footnote{The argument brakes down if one allows arbitrary
  $x\in\rn$ as replacing in the definition of $A_P$ the union over
  $x\in P$ by a union over $x\in\rn$ would yield $A_P(s)=0$.} 
To obtain the statement for arbitrary $x,y\in\rn$ we need to require that 
$d\varphi$ is \sP. 
\begin{thm}\label{thm-back}
Let $P$ be a Delone set of \flc. 
There exists $1>\epsilon>0$
such that for all differentiable functions
$\varphi:\rn\to\rn$ whose differential satisfies
$\|d(\varphi-\id)\|_\infty<\epsilon$ and is \sP\
holds:
$\forall r>0$ $\exists r'>0$ $\forall x,y\in \rn$
\begin{equation}
B_{r'}[\varphi(P)-\varphi(x)]=B_{r'}[\varphi(P)-\varphi(y)]\Longrightarrow
B_{r}[P-x]=B_{r}[P-y].
\end{equation}
\end{thm}
\begin{proof}
Given $P$ and $\epsilon$ according to Corollary~\ref{cor-back} let 
$\varphi:\rn\to\rn$ be a $C^1$-function
whose differential $d\varphi$ is \eP\ with range $R$ and satisfies
$\|d(\varphi-\id)\|_\infty<\epsilon$. Let $r\geq\max\{\lambda
r_{max},R\}$ where
$r_{max}=\inf\{r|\forall x\in\rn:B_r[P-x]\neq\emptyset\}$,
$\lambda=(1-\epsilon)^{-1}$. Let $r'=r'(r+\lambda r_{max})$ according
to Corollary~\ref{cor-back}. 

Suppose $B_{r'+\lambda r_{max}}[\varphi(P)-\varphi(x)]=B_{r'+\lambda
  r_{max}}[\varphi(P)-\varphi(y)]$. Choose $h\in
\varphi(P)-\varphi(x)$ such that $|h|\leq \lambda r_{max}$.
Let $x'=\varphi^{-1}(\varphi(x)+h)$ and
$y'=\varphi^{-1}(\varphi(y)+h)$. Then $\varphi(x'),\varphi(y')\in
\varphi(P)$ and $B_{r'}[\varphi(P)-\varphi(x')]=
B_{r'}[\varphi(P)-\varphi(y')]$. From Corollary~\ref{cor-back} we
conclude that $B_{r+\lambda r_{max}}[P-x']=B_{r+\lambda
  r_{max}}[P-y']$. Since $d\varphi$ is \eP\ with range $r$ the latter
implies that $\varphi(x)-\varphi(x')=\varphi(y'+x-x')-\varphi(y')$
which, in turn, is equivalent to  $\varphi(y)=\varphi(y'+x-x')$.
Since $\varphi$ is invertible we get $y-y'=x-x'$ and hence
$B_{r}[P-x]=B_{r}[P-y]$.
\end{proof}

\begin{cor}  \label{cor-back1}
Let $P$ be a Delone set of \flc. There exists $1>\epsilon>0$
such that for all differentiable functions
$\varphi:\rn\to\rn$ whose differential satisfies
$\|d(\varphi-\id)\|_\infty<\epsilon$ and is \sP\
we also have that $d\varphi^{-1}$ is
strongly $\varphi(P)$-equivariant.   
\end{cor}
\begin{proof}
Given $P$ there exists $\epsilon$ according to the last theorem.
Suppose $\|d(\varphi-\id)\|_\infty<\epsilon$ and that $d\varphi$ is
\eP\ with range $r$. Let $r'$ correspond to the one needed in the last theorem.
Let $B_{r'}(\varphi(P)-\varphi(x))=B_{r'}(\varphi(P)-\varphi(y))$.
By the last theorem this implies  $B_{r}(P-x)=B_{r}[P-y]$. 
Hence $d\varphi(x)=d\varphi(y)$. 
Since 
$$ ({d\varphi^{-1}})(\varphi(x))=(d\varphi(x))^{-1}$$
(the inverse in the matrix sense)
we see that $d\varphi^{-1}$ is \eP\ with range $r'$.
\end{proof}
We can now answer the question what it means if two closed forms differ by an
exact form, at least for forms which are close to $d\id$.  
\begin{cor} \label{cor-back2}
Let $P$ be a Delone set of \flc. 
There exists $1>\epsilon>0$
such that if  $d\varphi$ and $d\psi$ are \sP,
$\|d(\varphi-\id)\|_\infty<\epsilon$, 
and furthermore $\varphi-\psi$ \wP\ (or even \sP)
then $(\psi-\varphi)\circ\varphi^{-1}$ is weakly 
$\varphi(P)$-equivariant (or even strongly $\varphi(P)$-equivariant).
In particular, $\psi\circ\varphi^{-1}$ induces a topological semi-conjugacy
between $(\Omega_{\varphi(P)},\rn)$ and $(\Omega_{\psi(P)},\rn)$ (or even
$\psi(P)$ is locally derivable from $\varphi(P)$).
\end{cor}
\begin{proof} Suppose $\varphi-\psi$
  is \eP\ with range $r$. By
  Theorem~\ref{thm-back} we find $r'$ such that
  $B_{r'}[\varphi(P)-x]=  B_{r'}[\varphi(P)-y]$ implies 
$B_{r}[P-\varphi^{-1}(x)]= B_{r}[P-\varphi^{-1}(x)]$. 
Since $\eta(x)=(\varphi-\psi)(\varphi^{-1}(x))$ the latter implies
$\eta(x)=\eta(y)$. Hence $\eta$ is \eP\ with range $r'$ and so
the result follows directly from
Lemma~\ref{lem1}
(the $\varphi$, $g$ and $P$ there correspond to 
$\id-\eta$, $\id$ and $\varphi(P)$
here.)
\end{proof}

\subsection{Locally deriving back}
When are $\varphi(P)$ and $P$ mutually locally derivable?
Suppose that $\varphi-\id$ is \sP\ so that in particular $\varphi(P)$
is locally derivable from $P$. The last corollary gives an answer to
the question under which
condition $P$ can be locally derived back from $\varphi(P)$, namely if
$\|d(\varphi-\id)\|_\infty<\epsilon$ for some $\epsilon$ whose size
could be estimated using the details of Corollary~\ref{cor-back}.

\subsection{Cohomological interpretation}
We interprete our results in cohomological terms.

Let $B_{\|\cdot\|_\infty}(\alpha,\epsilon)$ be the open
$\epsilon$-ball around $\alpha\in Z^1_{s\dash P}(\rn,\rn)$ w.r.t.\ the
$\|\cdot\|_\infty$-norm. 
A good measure for the size of the deformation defined by $\varphi$ is
given by
$\inf\{\|d\varphi+d\eta\|_\infty: d\eta\in N^1_{s\dash  P}(\rn,\rn)\}$.
It is scale invariant and corresponds to the maximal relative change
of distance between points resulting from the deformation.

Theorem~\ref{thm-closed} and
Corollaries~\ref{cor-back1}, \ref{cor-back2} allow for the following
interpretation. 
For all Delone sets of \flc\ there
exists an $\epsilon$ such that the elements of $Z^1_{s\dash
  P}(\rn,\rn)\cap B_{\|\cdot\|_\infty}(d\id,\epsilon)$ define
invertible deformations of $P$. If moreover two such elements differ
by an element of $B^1_{w\dash P}(\rn,\rn)$ (or even $B^1_{s\dash
  P}(\rn,\rn)$) then they define pointed topological conjugate (or even
mutually locally derivable)
deformations. 

Theorems~\ref{thm-exact}, \ref{thm-exact-converse}
say that the elements near the class of $d\id$ in the mixed cohomology group
$ Z^1_{s\dash P}(\rn,\rn) / B^1_{w\dash P}(\rn,\rn) \cap Z^1_{s\dash
  P}(\rn,\rn)$
parametrise small deformations modulo bounded deformations which are 
in the same pointed conjugacy class. Here we say that a deformation
defined by a map $\varphi$ is bounded if $\varphi-\id$ is bounded, a
condition which implies that the deformation is homotopic to the
original pattern in the Hausdorff metric. Note that a bounded
deformation has finite size but the converse need not be true.
 
On the level of pattern spaces and their associated dynamical systems
the following picture has emerged (Theorem~\ref{thm-closed-converse}). 
An element $d\varphi$ near $d\id$ 
of $ Z^1_{s\dash P}(\rn,\rn)$ defines a homeomorphism
from the hull of $P$ to that of the deformation it defines 
which restricts to a
homeomorphism between the canonical transversals and maps $P$ to its
deformation. 
If $d(\varphi-\id)$ lies even in 
$ B^1_{w\dash P}(\rn,\rn) \cap Z^1_{s\dash P}(\rn,\rn)$
then it defines a second homeomorphism, this time a pointed
conjugacy between the continuous dynamical systems. 
The two homeomorphisms do in general not coincide.

\end{document}